\pdfoutput=1
\RequirePackage{ifpdf}
\ifpdf 
\documentclass[pdftex]{sigma}
\else
\documentclass{sigma}
\fi

\usepackage{cancel}

\numberwithin{equation}{section}

\newtheorem{Theorem}{Theorem}[section]
\newtheorem*{Theorem*}{Theorem}

\newtheorem{Lemma}[Theorem]{Lemma}

 { \theoremstyle{definition}

\newtheorem{Note}[Theorem]{Note}

\newtheorem{Remark}[Theorem]{Remark} }

\def\harr#1#2{\ \smash{\mathop{\hbox to .3in{\rightarrowfill}}\limits^{\scriptstyle#1}_{\scriptstyle#2}}\ }

\def\bra#1#2{\langle #1, #2\rangle}

\def\ss{\subset}

\def\arr{\longrightarrow}

\def\rn{\bbr^n}

\def\and{\qquad {\rm and} \qquad}
\def\arr{\longrightarrow}

\def\bbr{{\mathbb R}}\def\bbh{{\mathbb H}}
\def\bbc{{\mathbb C}}

\def\bbz{{\mathbb Z}}
\def\bbp{{\mathbb P}}

\def\l{\lambda}

\def\D{\Delta}

 \def\spin{{\rm Spin}}
 \def\spinh{{\rm Spin}^h}
 \def\cln{{\rm C}\ell_n}
 \def\clnh{{\rm C}\ell_{n, \bbh}}

\begin{document}
\allowdisplaybreaks

\renewcommand{\thefootnote}{}

\newcommand{\arXivNumber}{2301.09683}

\renewcommand{\PaperNumber}{012}

\FirstPageHeading

\ShortArticleName{Spin$^h$ Manifolds}

\ArticleName{Spin$\boldsymbol{{}^h}$ Manifolds\footnote{This paper is a~contribution to the Special Issue on Differential Geometry Inspired by Mathematical Physics in honor of Jean--Pierre Bourguignon for his 75th birthday. The~full collection is available at \href{https://www.emis.de/journals/SIGMA/Bourguignon.html}{https://www.emis.de/journals/SIGMA/Bourguignon.html}}}

\Author{H.~Blaine LAWSON Jr.}

\AuthorNameForHeading{H.B.~Lawson~Jr.}

\Address{Stony Brook University, Stony Brook NY, USA}
\Email{\href{mailto:blaine@math.sunysb.edu}{blaine@math.sunysb.edu}}
\URLaddress{\url{https://www.math.stonybrook.edu/~blaine/}}

\ArticleDates{Received January 25, 2023, in final form March 06, 2023; Published online March 19, 2023}

\Abstract{The concept of a ${\rm Spin}^h$-manifold, which is a cousin of Spin- and ${\rm Spin}^c$-manifolds, has been at the center of much research in recent years. This article discusses some of the highlights of this story.}


\Classification{53C27; 55P99}

\renewcommand{\thefootnote}{\arabic{footnote}}
\setcounter{footnote}{0}

It is a great pleasure for me to contribute a paper to this volume dedicated to Jean-Pierre Bourguignon
on his seventy-fifth birthday. Jean-Pierre has been a mathematical colleague and a dear friend
for most of my mathematical life, and throughout that time I have been in awe of his abilities to support
mathematics, and science in general, in so many ways.

Our various interests in mathematics have been close over the years, from manifolds of
negative curvature, minimal varieties, Yang--Mills fields, spin manifolds and Dirac operators.
I want to say that the best result in our collaboration on Yang--Mills \cite{Lawson:11, Lawson:10}
 was entirely Jean-Pierre's.

 In light of
his great interest in spin geometry, I thought it would be appropriate to recount
 some recent work in that area by various young people at Stony Brook.

 Let's begin by recalling that for $n\geq 3$ the group ${\rm Spin}(n)$ is the universal covering
 group
\[
\pi \colon \ {\rm Spin}(n) \to {\rm SO}(n)
\]
of the special orthogonal group. An oriented Riemannian manifold $X$ of dimension $n$ is a {\it Spin-manifold}
if there exists a principal ${\rm Spin}(n)$-bundle $P_{\rm Spin} \to X$ on $X$, and a ${\rm Spin}(n)$-equivariant bundle map
\[
 P_{\rm Spin} \to P_{\rm SO}(X),
\]
where $ P_{\rm SO}(X)$ is the bundle of positively oriented, orthonormal tangent frames on $X$.
Each such ``square root'' of the principal bundle $P_{\rm SO}(X)$, up to equivalence,
is called a {\it Spin-structure} on~$X$. The existence of such structures depends only
on the topology of~$X$. That is, Spin-structures exist for any Riemannian metric on $X$
if and only if the first two Stiefel--Whitney classes satisfy $w_1(X)= w_2(X)= 0$.

Now the notion of a Spin-manifold played an important role in the establishment of the
Atiyah--Singer index theorem. It was a result of Atiyah and Hirzebruch \cite{Lawson:6} that every compact
Spin-manifold has the property that its $\hat{\mathbb A}$-genus is an integer. This characteristic invariant
is generally a rational number. (For example, for the complex projective plane $\bbp^2(\bbc)$, it is $-1/8$.)
Atiyah thought that this might be a consequence of the fact that $\hat{\mathbb A}(X)$
was the index of an elliptic operator on $X$ that could not be defined on non-Spin manifolds.
Indeed Singer found such an operator, and this lit the path to the theorem and its proof.

These operators were obtained by defining vector bundles on $X$ as associated bundles
\[
{\mathbf V}_\l \equiv P_{\rm Spin} \times_\lambda V \to X,
\]
 where $\l \colon {\rm Spin}(n) \to {\rm SO}(V)$ is a finite-dimensional representation of ${\rm Spin}(n)$ which
 is {\it not} pulled back from a representation of SO$(n)$. Such representations can be obtained from
 representations of the Clifford algebra ${\rm C}\ell(\rn, \bra \cdot \cdot)$.

We recall \cite{Lawson:17} that a Riemannian manifold $(X, \bra\cdot \cdot)$ has a natural bundle of Clifford algebras~${\rm C}\ell(X)$ whose fibre at $x\in X$ is the Clifford algebra~${\rm C}\ell(T_xX, \bra \cdot \cdot)$.
For Clifford algebra representations $\l$, the associated bundle ${\mathbf V}_\l $ is a bundle of
modules for the bundle of Clifford algebras~${\rm C}\ell(X)$. This fact allows one to define differential operators
via the Dirac construction.

It turns out that the notion of a Spin-manifold has analogues over the complex numbers and the quaternions.
In the complex case they are called ${\rm Spin}^c$-manifolds, and they have been of central importance for quite some time.
Those in the quaternion case I will call ${\rm Spin}^h$-manifolds, and they will be the center of discussion for the rest of the paper.

 To begin we shall consider just the groups themselves (or, equivalently, suppose the manifold to be just a point).
 Let's start by writing ${\rm Spin}(n)$ trivially as ${\rm Spin}(n)\times_{\bbz_2} \bbz_2 = {\rm Spin}(n)\times_{\bbz_2} {\rm O}(1)$
 where $\bbz_2$ acts diagonally. Then we can define (with $\bbz_2$ always acting diagonally)
 \begin{gather*}
{\rm Spin}^r (n) = {\rm Spin} (n) \times_{\bbz_2} {\rm O}(1), \\
 {\rm Spin}^c (n) = {\rm Spin} (n) \times_{\bbz_2} {\rm U}(1), \\
{\rm Spin}^h (n) = {\rm Spin} (n) \times_{\bbz_2} {\rm Sp}(1),
 \end{gather*}
 the {\it real, complex and quaternionic Spin groups}.
 Dividing further by the group $\bbz_2\times\bbz_2$ (where the original $\bbz_2$ is the diagonal) gives us surjective homomorphisms
 \begin{gather*}
 {\rm Spin}^r (n) \arr {\rm SO} (n), \\
 {\rm Spin}^c (n) \arr {\rm SO} (n) \times {\rm U}(1), \\
 {\rm Spin}^h (n) \arr {\rm SO} (n) \times {\rm SO}(3).
 \end{gather*}
Now a {\it ${\rm Spin}^c$-manifold} is an oriented Riemannian $n$-manifold, equipped with a principal
 ${\rm Spin}^c(n)$-bundle $P_{{\rm Spin}^c}$, and a ${\rm Spin}^c(n)$-equivariant bundle map
\[
 P_{{\rm Spin}^c} \to P_{{\rm SO}}(X) \times P_{{\rm U}(1)}(L), \qquad \text{ (fibre product)}
\]
 where $L$ is a complex line bundle on $X$, called the {\it canonical bundle} of the structure, and $P_{{\rm SO}}(X)$ is as above.
 Such a structure exists with canonical bundle $L$ if and only if $w_2(X)=\rho(c_1(L))=0$
 where $\rho$ is mod-2 reduction (cf.\ \cite[p.~391]{Lawson:17}).
 Now complex line bundles on $X$ are in 1-to-1 correspondence with elements of $H^2(X,\bbz)$.
 Thus a manifold carries a ${\rm Spin}^c$-structure if and only if $w_2(X)$ is the mod-2 reduction of an integer class,
 which is equivalent to the third integral Stiefel--Whitney class satisfying $W_3(X)=0$. Note that given one ${\rm Spin}^c$-structure on $X$,
 we get all the rest by tensoring $L$ with complex line bundles having $w_2=0$.
 We point out that every almost complex manifold $X$ has a canonical ${\rm Spin}^c$-structure where $L$ is the anticanonical bundle~of $X$.

From this point of view ${\rm Spin}^h$-manifolds are defined analogously.
 A {\it ${\rm Spin}^h$-manifold} is an oriented Riemannian $n$-manifold, equipped with a principal
 ${\rm Spin}^h(n)$ bundle $P_{{\rm Spin}^h}$, and a ${\rm Spin}^h(n)$-equivariant bundle map
\[
 P_{{\rm Spin}^h} \to P_{{\rm SO}}(X) \times P_{{\rm SO}(3)}(E),
\]
 where $E$ is an oriented Riemannian 3-plane bundle on $X$, called
 the {\it canonical bundle} of the structure. Such a structure exists if and only if one
 can find such an $E$ with $w_2(X)+w_2(E)=0$.

Now for $\spinh$-manifolds there is an analogy with the fact that an almost complex manifold is ${\rm Spin}^c$.
 A manifold of dimension $4n$ is {\it almost quaternionic} if its structure group can be reduced to
 ${\rm Sp}(n)\cdot {\rm Sp}(1) = {\rm Sp}(n)\times_{\bbz_2} {\rm Sp}(1)$. Every such manifold
 has a canonical ${\rm Spin}^h$-structure. In dimensions 8$n$, the manifold is in fact spin.
 In dimensions~4 (mod~8) the canonical bundle $E$ is a 3-dimensional subbundle of ${\rm End}(TX)$
 which locally has an orthonormal basis $I$, $J$, $K$ with $I^2=J^2=K^2= IJK=-{\rm Id}$. See \cite[Remark~2.2]{Lawson:1}
 for a detailed discussion.

 Note there are natural homomorphisms ${\rm Spin}(n) \to {\rm Spin}^c(n) \to {\rm Spin}^h(n)$.
 This shows that any Spin-manifold has a ${\rm Spin}^c$-structure by choosing the canonical bundle to be
 trivial, and any ${\rm Spin}^c$-manifold has a ${\rm Spin}^h$-structure with canonical bundle $L\oplus \bbr$.
 Furthermore, the notion of a~${\rm Spin}^c$-manifold is definitely more general than Spin,
 since $\bbp^2(\bbc)$ is complex and therefore ${\rm Spin}^c$, but not Spin,
 as we have discussed.

In a similar vein there are manifolds which are not ${\rm Spin}^c$ but are ${\rm Spin}^h$. A basic example is
the Wu manifold $X={\rm SU}(3)/{\rm SO}(3)$ which is not ${\rm Spin}^c$ ($W_3\neq 0$), but is ${\rm Spin}^h$.
This was first noticed by Xuan Chen in his thesis \cite[p.~26ff]{Lawson:12}.
The general question of understanding manifolds with this property is treated in great generality in the work of
Albanese and Milivojevi\'c below.

Now just like Spin- and ${\rm Spin}^c$-manifolds, ${\rm Spin}^h$-manifolds also carry Dirac operators.
 Let $\cln$ denote the Clifford algebra of $\rn$ with the standard inner product, and set
\[
 \clnh \equiv \cln \otimes_\bbr \bbh.
\]
 Then there is an embedding
\[
 \spinh(n) \ss \clnh
\]
 given as the projection of $\spin(n) \times {\rm Sp}(1) \ss \cln\times \bbh$.
 Suppose now that $X$ is a $\spinh$-manifold. Then given any representation
 $\rho\colon \clnh \to {\rm End}(V)$ on a finite-dimensional vector space $V$, we restrict to $\spinh$
 and form the associated bundle
\[
 {\bf V}_\rho = P_{\spinh} \times_{\rho} V.
\]
 If the representation $\rho$ is irreducible, the bundle $ {\bf V}_\rho $ is called {\it fundamental}.
 There is always a~Dirac operator defined on sections of ${\bf V}_\rho$.
 In certain dimensions the restriction of an irreducible~$\rho$ to the even part of~$\clnh$
 splits into two irreducible representations, and the Dirac operator interchanges the
 sections of the corresponding bundles. For example, this happens when $n\equiv 0$ $({\rm mod} \ 8)$.
 For this question, and other as well, there is a periodicity of order~8 in dimensions, as in the real case.
 However one must be careful. An irreducible real representation of~$\cln$, when tensored with $\bbh$, may not be an
 irreducible real representation of~$\clnh$.

 There are $\clnh$-linear Dirac operators with indices in ${\rm KSp}_*$ which has been worked out
 by Jiahao Hu~\cite{Lawson:16} and is discussed below.

There is an interesting story surrounding Spin-manifolds and their cousins.
 It was a ``folk'' theorem for some years that an oriented Riemannian manifold $X$
has a Spin-structure if and only if there exists a vector bundle of real ${\rm C}\ell(X)$-modules
$\cancel{S} \to X$ which are irreducible at each point. That such a bundle
$\cancel{S}$ exists on a Spin-manifold was clear, and I asked my graduate student Xuan Chen, some years ago,
 to find a good geometric proof of the converse.

Xuan first showed that if one replaces Spin-structures
with ${\rm Spin}^c$-structures, and the existence of a real locally irreducible
${\rm C}\ell(X)$-module with the existence of a complex one, then
 the converse is always true.
Furthermore, the proof of this result has a very geometric flavor.
Xuan also found a proof in the real case whenever the dimension of~$X$ is
congruent to 0, 6 or 7 (mod 8). However, to my great surprise he showed that
{\it in all other dimensions $> 8$ the converse is false!}

Xuan then came up with the notion of a ${\rm Spin}^h$-structure on a manifold,
as explained above. He showed that the canonical bundles $E$ are characterized by the condition
 on the Stiefel--Whitney classes that: $w_2(TX) = w_2(E)$.
He then showed that ${\rm Spin}^h$-structures can exist on manifolds
which have no ${\rm Spin}^c$-structure, and in
 dimensions 2, 3, or 4 (mod 8), a ${\rm Spin}^h$-structure
implies the existence of a real irreducible ${\rm C}\ell(X)$-module.

I found this idea of a quaternionic analogue of a ${\rm Spin}^c$-manifold enticing.
 Michael Albanese and Aleksandar Milivojevi\'c, students at Stony Brook at the time, showed that
the primary obstruction to the existence of a ${\rm Spin}^h$-structure is the fifth integral
Stiefel--Whitney class $W_5$. They then showed that:
\begin{enumerate}\itemsep=0pt
\item[] {\it Every compact oriented manifold of dimension $\leq 7$ admits a ${\rm Spin}^h$-structure.}
\item[] {\it However, in every dimension $\geq 8$,
there are infinitely many homotopy types of compact
simply connected manifolds which do not admit ${\rm Spin}^h$-structures.}
\end{enumerate}
They have similar results for non-compact manifolds:
\begin{enumerate}\itemsep=0pt
\item[] {\it Every oriented $5$-manifold is ${\rm Spin}^h$.}
\end{enumerate}
\emph{Note:} Every oriented 4-manifold is ${\rm Spin}^c$ \cite{Lawson:15} (see also \cite{Lawson:19, Lawson:20}, for compact manifolds, and the case of
non-compact manifolds was done
in~\cite{Lawson:25}), and every oriented 3-manifold is Spin:
\begin{enumerate}\itemsep=0pt
\item[] {\it An oriented manifold of dimension 6 or 7 is ${\rm Spin}^h$ if and only if
the primary obstruction $W_5$ vanishes.}
\end{enumerate}
See \cite{Lawson:1,Lawson:2} for these results, and \cite{Lawson:2} for further discussion of the non-compact case.

It turned out that the notion of ${\rm Spin}^h$-manifolds had been found previously
(but the results mentioned here are new). See, for example, \cite{Lawson:9, Lawson:21}.
Recent math papers include \cite{Lawson:1,Lawson:2, Lawson:12,Lawson:13,Lawson:18}.

Now in fact, it turned out that {\it physicists
had also found ${\rm Spin}^h$-manifolds}, and were interested in them for physical reasons.
See, for example, \cite{Lawson:13,Lawson:23}, and in dimension 4 \cite{Lawson:7,Lawson:8,Lawson:14,Lawson:26}.
In \cite{Lawson:22}, a ${\rm Spin}^h$-structure on a 4-manifold was used in constructing a quaternionic version of Seiberg--Witten theory.

In fact the work in~\cite{Lawson:1} goes much further than my ${\rm Spin}^h$ discussion. They define a notion of a ${\rm Spin}^k$-manifold
for all integers $k>0$ where $k=1,2,3$ correspond to the three cases above. For $k\geq 3$ the group is defined by
\[
{\rm Spin}^k(n) \equiv {\rm Spin}(n) \times_{\bbz_2} {\rm Spin}(k)
\]
with $\bbz_2$ acting diagonally. They prove that:
\begin{enumerate}\itemsep=0pt
\item[] {\it For every $k$ there exists a compact simply connected smooth manifold
which does not admit a ${\rm Spin}^k$-structure.}
\end{enumerate}
Their results on ${\rm Spin}^k$ manifolds are non-trivial and basic for their results in the ${\rm Spin}^h$ case.
They also have the corollary that:
\begin{enumerate}\itemsep=0pt
\item[] {\it There is no integer $k$ such that
 every manifold embeds into a
Spin-manifold with codimension $k$.}
\end{enumerate}
Note that:
\begin{enumerate}\itemsep=0pt
\item[] {\it Every manifold $X$ embeds into a orientable manifold with codimension $1$}
\end{enumerate}
by embedding $X$ as the zero-section of its orientation line bundle. Asking for the spin analogue
of this was part of the motivation for examining ${\rm Spin}^k$-manifolds, and the answer they found was
striking.

Recently ${\rm Spin}^h$-manifolds played a role in the work of Freed and Hopkins~\cite{Lawson:13} which discusses
quantum field theories, Thom spectra, and much else. In their paper, ${\rm Spin}^h$-manifolds,
Thom classes of ${\rm Spin}^h$-bundles, ${\rm Spin}^h$-cobordism,
 its natural transformation to symplectic ${\rm K}$-theory, etc.\ are all defined, but as special cases
of a more general discussion. The ${\rm Spin}^h$-case appears, for example, in case~4 of Table~9.35.
Certainly this shows that ${\rm Spin}^h$-theory is relevant in modern quantum field theories.

Very recently, with the work of Jiahao Hu in his Stony Brook thesis \cite{Lawson:16},
{\it ${\rm Spin}^h$-manifolds have become a basic tool in topology.}

 This resulted from his addressing the following problem:
 {\it Find enough invariants to know that two real vector bundles over a space $X$ are stably isomorphic.}
 Of course the bundles must have the same characteristic classes, but that is not enough.
 (For example, there exists a non-stably trivial real vector bundle over $S^9$ with vanishing
Stiefel--Whitney and Pontryagin classes.)
 Jiahao Hu has solved this problem and done it in a very geometric way, using ``integration over cycles''.

 Now stable equivalence classes of real bundles on a reasonable space $X$ form the real ${\rm K}$-ring ${\rm KO}(X)$
 which expands to a generalized cohomology theory ${\rm KO}^*(X)$. One of Hu's important contributions
 was showing that:
\begin{enumerate}\itemsep=0pt
\item[] {\it Classes in a generalized cohomology theory are determined by
their periods over cycles for the ``Anderson dual theory''.}
\end{enumerate}

Now the Anderson dual of ${\rm KO}$-theory turns out to be symplectic ${\rm K}$-theory
${\rm KSp}^*$, which is defined like ${\rm KO}^*$ but by using quaternionic vector bundles in place of real ones~\cite{Lawson:3}.

 This led Hu to a study of symplectic ${\rm K}$-theory and ${\rm Spin}^h$-cobordism with many new results.
 One of his basic theorems is that:
\begin{enumerate}\itemsep=0pt
\item[] {\it The Atiyah--Bott--Shapiro isomorphism~{\rm \cite{Lawson:4}} holds
for quaternionic Clifford modules and symplectic ${\rm K}$-theory.}
\end{enumerate}

 Jiahao then defined a {\it Thom class} for ${\rm Spin}^h$-vector bundles in symplectic K-theory,
 using $(\D^+-\D^-)\otimes_\bbr \bbh$ where $\D^\pm$
 are the real irreducible spinor bundles for
 ${\rm Spin}(8n)$ and $\bbh$ is the canonical representation of~${\rm Sp}(1)$.
 He constructed the Thom spectrum ${\rm MSpin}^h$ for ${\rm Spin}^h$-cobordism, and defined a map
 \[
 {\rm ind}^h \colon \ {\rm MSpin}^h \arr {\rm KSp},
\]
 which allowed him to define the cycles needed for real ${\rm K}$-theory.
 (As he found out later, and as mentioned above, these things had already been done
 in~\cite{Lawson:13}.)

 \def\sshh{{\mathbb S}}

 His next important step was to show that:
\begin{enumerate}\itemsep=0pt
\item[] {\it There exist Dirac operators on ${\rm Spin}^h$-manifolds
 with indices in ${\rm KSp}_*({\rm pt})$.}
\end{enumerate}

 To do this he followed Atiyah--Singer's index theorem for elliptic operators on ${\rm C}\ell_k$-bundles,
 where the algebra ${\rm C}\ell_k$ acts universally on the bundle~\cite{Lawson:5}. (One could also see \cite[p.~139ff]{Lawson:17}.)
 In this context, for the Dirac operator on a Spin $n$-manifold $X$, one defines a
 ${\rm C}\ell_n$-bundle
 \[
 \sshh(X) = P_{{\rm Spin}}(X) \times_{\ell} {\rm C}\ell_n,
\]
 where $\ell \colon {\rm Spin}(n) \to {\rm End}( {\rm C}\ell_n)$ is given by left multiplication,
 and $ {\rm C}\ell_n$ acts on the bundle by right multiplication.
 The kernel is a real $\bbz_2$-graded $ {\rm C}\ell_n$-module which can be identified
 with ${\rm KO}_{n}$({\rm pt}) by the ABS isomorphism.

 On a ${\rm Spin}^h$-manifold $X$ one defines
\[
 \sshh_{\bbh}(X) = P_{{\rm Spin}^h}(X) \times_{\ell} ({\rm C}\ell_n \otimes_\bbr \bbh),
\]
 where $\ell\colon \spinh \equiv \spin(n) \times_{\bbz_2} {\rm Sp}(1) \arr {\rm End}(\cln \otimes_\bbr \bbh)$
 is given by left Clifford multiplication of Spin($n$) on $\cln$ and of Sp(1) on $\bbh$.
 This is a bundle of $\bbz_2$-graded $\cln \otimes_\bbr \bbh$-modules.
 The standard Dirac operator commutes with this action of $\cln \otimes_\bbr \bbh$,
 and so its analytic index lies in the Grothendieck group of $\bbz_2$-graded $(\cln \otimes_\bbr \bbh)$-modules,
 which is isomorphic to ${\rm KSp}_{n}$({\rm pt}) by the quaternionic version of the ABS isomorphism.
 One now needs the symplectic version of the Atiyah--Singer index theorem, which Jiahao provides in~\cite{Lawson:16}.

 From here Jiahao was able to prove the following:
\begin{enumerate}\itemsep=0pt
\item[] {\it
 ${\rm Spin}^h$-manifolds provide enough cycles for symplectic ${\rm K}$-theory
to distinguish real vector bundles up to stable equivalence. }
\end{enumerate}

 By a {\it Spin-cycle} on a space $Z$ (a finite CW-complex)
 we mean a map $f\colon X\to Z$ where $X$ is a compact Spin-manifold. For a real vector bundle
 $E\to Z$ we define its period over the Spin-cycle $f$ to be
 $\bra f E \equiv$ the ${\rm KO}$-index
 of the Dirac operator on $X$ twisted by $f^*E$,
 which takes values in $\bbz$ or $\bbz_2$ depending on dimension.
 A {\it ${\rm Spin}^h$-cycle} and the period of a real vector bundle over a ${\rm Spin}^h$-cycle
 are defined analogously using ${\rm Spin}^h$-manifolds and the Dirac operator on
 ${\rm Spin}^h$-manifolds.

 We think of these as integer invariants. There are also torsion invariants defined in
 much the same way.
The assertion above says that the periods of a real vector bundle $E$ over ${\rm Spin}^h$-cycles
 and torsion ${\rm Spin}^h$-cycles determine the stable isomorphism class of~$E$.

The proof uses the fact that Jiahao's quaternionic version of the Atiyah--Bott--Shapiro isomorphism is
 equivariant with respect to the real ABS isomorphism.
 This implies that ${\rm ind}^h$ above is equivariant with
 respect to the Atiyah--Bott--Shapiro orientation ${\rm ind}\colon {\rm MSpin} \to {\rm KO}$.
 Since ind$^h\big(\bbp^1(\bbh)\big)$ generates ${\rm KSp}_4({\rm pt})\cong \bbz$,
 one can transfer periods of real bundles defined over Spin-cycles $f\colon X\to Z$
 to periods over ${\rm Spin}^h$-cycles $f\times f_0\colon X \times \bbp^1(\bbh) \to Z \times{\rm pt}$.

 This work gives an answer to the question posed above. The history of this question
 goes back to Dennis Sullivan's study of Hauptvermutung for manifolds~\cite{Lawson:24}. For that purpose,
 he analyzed ordinary cohomology, complex ${\rm K}$-theory and real ${\rm K}$-theory at odd primes.
 His arguments relied on knowing that these theories are Anderson self-dual and that,
 for every $p$,
 the associated theories with $\bbz_p$-coefficients have a cup product.
 Both of these fail for real ${\rm K}$-theory at prime~2.
 Jiahao overcame this difficulty.
 His result at the prime~2 is really surprising.

 It is common in topology that prime 2 is more difficult than odd primes. Jiahao formulated and proved a general statement about determining cohomology classes in a generalized cohomology theory by their periods over cycles for its Anderson dual theory.
 In the real ${\rm K}$-theory case, he needed cycles for symplectic ${\rm K}$-theory to define invariants taking values in $\bbz$ and $\bbz_2$,
 which forced him to consider ${\rm Spin}^h$-manifolds and to define the KSp-index of Dirac operators on ${\rm Spin}^h$-manifolds.

 In another direction, Jiahao has also computed the $\bbz_2$-cohomology of the classifying space ${\rm BSpin}^h$. It
 turns out to be very similar to the $\bbz_2$-cohomology of BSpin and ${\rm BSpin}^c$.
 He found, interestingly, that
 \begin{gather*}
 H^*({\rm BSpin};\bbz_2) = H^*({\rm BSO};\bbz_2)/\big(\nu_2, Sq^1 \nu_{2^r} , r \geq 1\big) \\
 H^*\big({\rm BSpin^c};\bbz_2\big) = H^*({\rm BSO};\bbz_2)/\big(Sq^1 \nu_{2^r} , r \geq 1\big) \\
 H^*\big({\rm BSpin^h};\bbz_2\big) = H^*({\rm BSO};\bbz_2)/\big(Sq^1 \nu_{2^r} , r \geq 2\big),
 \end{gather*}
 where $\nu_k$ is the $k^{\rm th}$ Wu class.
 He has also computed the first several ${\rm Spin}^h$ cobordism groups.

I hope I have been able to attract your interest to ${\rm Spin}^h$-manifolds. I think they play an important
 role in understanding certain aspects of geometry.

\subsection*{Acknowledgements}
I want to thank Michael Albanese, Aleksandar Milovojevi\'c and Jiahao Hu for their careful reading of this
 manuscript and for the many serious improvements they have suggested.
 I would like to thank the Simons Foundation for support during the writing of this paper.

\pdfbookmark[1]{References}{ref}
\LastPageEnding

\end{document}